\documentclass{gtart_h}

\def\ifplaintex{\expandafter\ifx\csname documentclass\endcsname\relax}

\ifplaintex 
\else
\fi

\def\gt{{\mathsurround=0pt\it $\cal G\mskip-2mu$eometry \&\ 
$\cal T\!\!$opology}}        

\def\gtm{{\mathsurround=0pt\it $\cal G\mskip-2mu$eometry \&\ 
$\cal T\!\!$opology $\cal M\mskip-1mu$onographs}}    

\def\gtp{{\mathsurround=0pt\it $\cal G\mskip-2mu$eometry \&\ 
$\cal T\!\!$opology $\cal P\!$ublications}}  

\def\recd{{\small Received:\qua\receiveddate\ifx\reviseddate\relax
\else\qquad Revised:\qua\reviseddate\fi\par}} 

\def\volumenumber#1{\def\thevolumenumber{#1}}
\def\volumeyear#1{\def\thevolumeyear{#1}}
\def\volumename#1{\def\thevolumename{#1}}
\def\papernumber#1{\def\thepapernumber{#1}}
\def\pagenumbers#1#2{\def\startpage{#1}\def\finishpage{#2}}
\def\published#1{\def\publishdate{#1}}
\def\received#1{\def\receiveddate{#1}}
\def\revised#1{\def\reviseddate{#1}}
\def\accepted#1{\def\accepteddate{#1}}

\long\def\asciiabstract#1{\long\def\theasciiabstract{#1}}
\def\asciikeywords#1{\def\theasciikeywords{#1}}

\let\\\par
\let\thevolumenumber\relax\let\thepapernumber\relax
\let\thevolumeyear\relax\let\startpage\relax
\let\finishpage\relax\let\publishdate\relax\let\receiveddate\relax
\let\reviseddate\relax\let\accepteddate\relax\let\theasciititle\relax
\let\theasciiauthors\relax
\let\theasciiabstract\relax\let\theasciikeywords\relax

\let\theerratum\relax\let\theasciiemail\relax
\let\theshortauthors\relax\let\theshorttitle\relax

\def\startpage{1}\def\finishpage{15}\def\thepapernumber{77}

\volumenumber{2}
\volumename{Proceedings of the Kirbyfest}
\volumeyear{1999}

\long\def\maketitlep{   

\count0=\startpage

\gtm\nl        
{\small Volume \thevolumenumber: \thevolumename\nl 
\ifx\theerratum\relax\else Erratum \erratumnumber\nl\fi
Pages \startpage--\finishpage\nl}

\vglue 0.1truein   

{\parskip=0pt\leftskip 0pt plus 1fil\def\\{\par\smallskip}{\ifplaintex\large
\else\Large\fi\bf\thetitle}\par\medskip}   
\vglue 0.05truein 

{\parskip=0pt\leftskip 0pt plus 1fil\def\\{\par}{\sc\theauthors}
\par\medskip}%
 
\vglue 0.03truein 

{\small\leftskip 25pt\rightskip 25pt{\bf Abstract}\stdspace\theabstract

{\bf AMS Classification}\stdspace\theprimaryclass
\ifx\thesecondaryclass\relax\else; \thesecondaryclass\fi\par
{\bf Keywords}\stdspace \thekeywords\par}\vglue 7pt

}   

\font\phead=cmsl9 scaled 950
\font=cmsl9 scaled 1050
\font\pnum=cmbx10 scaled 913
\font\lnum=cmbx10 
\font\pfoot=cmsl9 scaled 950
\font=cmsl9 scaled 1050
\ifplaintex
\headline{\vbox to 0pt{\vskip -4.5mm\line{\small\phead\ifnum
\count0=\startpage ISSN 1464-8997 (on line)
1464-8989 (printed) \hfill {\pnum\folio}\else\ifodd\count0\def\\{ }%
\ifx\theshorttitle\relax\thetitle\else\theshorttitle\fi\hfill{\pnum\folio}
\else\def\\{ and }{\pnum\folio}\hfill\ifx\theshortauthors\relax\theauthors
\else\theshortauthors\fi\fi\fi}\vss}}
\footline{\vbox to 0pt{\vglue 0mm\line{\small\pfoot\ifnum\count0=\startpage
Published \publishdate:\qua\copyright\ \gtp\hfill\else
\gtm, Volume \thevolumenumber\ (\thevolumeyear)\hfill\fi}\vss
}}
\else
\makeatletter
\def\@oddhead{{\small\lhead\ifnum\count0=\startpage ISSN 1464-8997 (on line)
1464-8989 (printed) \hfill {\lnum\number\count0}\else\ifodd\count0
\def\\{ }\ifx\theshorttitle\relax \thetitle \else\theshorttitle\fi\hfill
{\lnum\number\count0}\else\def\\{ and }{\lnum\number\count0}
\hfill\ifx\theshortauthors\relax 
\theauthors\else\theshortauthors\fi\fi\fi}}\def\@evenhead{@oddhead}
\def\@oddfoot{\small\lfoot\ifnum\count0=\startpage Published \publishdate:\qua\copyright\ \gtp\hfill\else
\gtm, Volume \thevolumenumber\ (\thevolumeyear)\hfill\fi}
\def\@evenfoot{@oddfoot}
\makeatother
\fi

\let\maketitlepage\maketitlep
\let\makeshorttitle\maketitlepage
\let\maketitle\maketitlepage

\newwrite\gtoutfile
\long\gdef\makeheadfile{  
{\def\\{, }\def\s{ }
\immediate\openout\gtoutfile head.xxx
\immediate\write\gtoutfile{Proxy-for: \ifx\theasciiauthors\relax
\theauthors\else\theasciiauthors\fi\s<\ifx\theasciiemail\relax\theemail\else\theasciiemail\fi>}
\immediate\write\gtoutfile{\noexpand\\}
\immediate\write\gtoutfile{Authors: \ifx\theasciiauthors\relax
\theauthors\else\theasciiauthors\fi}
{\def\\{ }\immediate\write\gtoutfile{Title: \ifx\theasciititle\relax
\thetitle\else\theasciititle\fi}}
\immediate\write\gtoutfile{Subj-class: GT or SG, GR etc}
\immediate\write\gtoutfile{MSC-class: \theprimaryclass\ifx\thesecondaryclass\relax\else, \thesecondaryclass\fi}
\immediate\write\gtoutfile{Journal-ref: Geom. Topol. Monogr. \thevolumenumber\s
(\thevolumeyear) \startpage-\finishpage}
\immediate\write\gtoutfile{Comments: Published by Geometry and Topology Monographs at}
\immediate\write\gtoutfile{\s\s\s  http://www.maths.warwick.ac.uk/gt/GTMon\thevolumenumber/paper\thepapernumber.abs.html}
\immediate\write\gtoutfile{\noexpand\\}
\immediate\write\gtoutfile{}
\ifx\theasciiabstract\relax
\immediate\write\gtoutfile{\theabstract}\else
\immediate\write\gtoutfile{\theasciiabstract}\fi
\immediate\write\gtoutfile{}
\immediate\write\gtoutfile{\noexpand\\}
\immediate\write\gtoutfile{}
\immediate\closeout\gtoutfile}}  

\def\maketitlepage{\maketitlep\makeheadfile}
\let\makeshorttitle\maketitlepage
\let\maketitle\maketitlepage

\volumenumber{7}
\volumename{Proceedings of the Casson Fest}
\volumeyear{2004}
\papernumber{12}
\pagenumbers{291}{309}
\received{16 June 2003}
\revised{1 November 2003}
\accepted{15 December 2003}
\published{20 September 2004}

\usepackage{amssymb,amsmath}

\theoremstyle{plain}

\newtheorem{proposition}{Proposition}[section]
\newtheorem{lemma}[proposition]{Lemma}
\newtheorem{corollary}[proposition]{Corollary}

\newtheorem{conjecture}{Conjecture}

\theoremstyle{definition}
\newtheorem{definition}[proposition]{Definition}

\newtheorem{question}{Question}

\theoremstyle{remark}

\newtheorem{remark}[proposition]{Remark}

\def\lbl{\label}

\def\BN{\mathbb N}
\def\BZ{\mathbb Z}

\def\BQ{\mathbb Q}

\def\BC{\mathbb C}

\def\BK{\mathbb K}

\def\A{\mathcal A}
\def\B{\mathcal B}

\def\F{\mathcal F}

\def\calI{\mathcal I}

\def\La{\Lambda}
\def\l{\lambda}

\def\s{\sigma}

\def\la{\langle}
\def\ra{\rangle}

\def\e{\epsilon}

\def\s{\sigma}

\def\sub{\subset}

\def\ch{\mathrm{ch}}

\def\loc{\mathrm{loc}}
\def\Hom{\mathrm{Hom}}

\def\longto{\longrightarrow}

\def\pt{\partial}
\def\fg{\mathfrak{g}}

\def\SL{\mathrm{SL}}
\def\ord{\mathrm{ord}}

\begin{document}


\title{On the characteristic and deformation varieties\\of a knot}

\author{Stavros Garoufalidis}
\address{School of Mathematics, Georgia Institute of 
Technology\\Atlanta, GA 30332-0160, USA}
\urladdr{http://www.math.gatech.edu/~stavros/}
\email{stavros@math.gatech.edu}
\primaryclass{57N10}
\secondaryclass{57M25}
\keywords{$q$--holonomic functions, $D$--modules, characteristic
variety, deformation variety, colored Jones function, multisums, 
hypergeometric functions, WZ algorithm.}

\asciikeywords{q-holonomic functions, D-modules, characteristic
variety, deformation variety, colored Jones function, multisums, 
hypergeometric functions, WZ algorithm.}

\begin{abstract}
The colored Jones function of a knot is a sequence of Laurent polynomials
in one variable, whose $n$th term is the Jones polynomial of the knot
colored with the 
$n$--dimensional irreducible representation of $\mathfrak{sl}_2$.
It was recently shown by TTQ Le and the author that the colored Jones
function of a knot is $q$--holonomic, ie, that it satisfies a nontrivial
linear recursion relation with appropriate coefficients. Using holonomicity, 
we introduce a geometric 
invariant of a knot: the characteristic variety, an affine 1--dimensional
variety in $\BC^2$.  We then compare it with the character variety of 
$\SL_2(\BC)$ representations, viewed from the boundary.
The comparison is stated as a conjecture which we verify (by a direct 
computation) in the case of the trefoil and figure eight knots. 

We also propose a geometric relation between
the peripheral subgroup  of the knot group, and basic operators
that act on the colored Jones function. We also define a noncommutative
version (the so-called noncommutative $A$--polynomial) of the characteristic 
variety of a knot.

Holonomicity works well for higher rank groups and goes beyond hyperbolic
geometry, as we explain in the last chapter.
\end{abstract}

\asciiabstract{%
The colored Jones function of a knot is a sequence of Laurent
polynomials in one variable, whose n-th term is the Jones polynomial
of the knot colored with the n-dimensional irreducible representation
of SL(2).  It was recently shown by TTQ Le and the author that the
colored Jones function of a knot is q-holonomic, ie, that it satisfies
a nontrivial linear recursion relation with appropriate coefficients.
Using holonomicity, we introduce a geometric invariant of a knot: the
characteristic variety, an affine 1-dimensional variety in C^2.  We
then compare it with the character variety of SL_2(C) representations,
viewed from the boundary.  The comparison is stated as a conjecture
which we verify (by a direct computation) in the case of the trefoil
and figure eight knots.

We also propose a geometric relation between
the peripheral subgroup  of the knot group, and basic operators
that act on the colored Jones function. We also define a noncommutative
version (the so-called noncommutative A-polynomial) of the characteristic 
variety of a knot.

Holonomicity works well for higher rank groups and goes beyond hyperbolic
geometry, as we explain in the last chapter.}

\maketitle

\cl{\small\it Dedicated to Andrew Casson on the occasion of his 60th
birthday}



\section{Introduction}
\lbl{sec.intro}

\subsection{The colored Jones function of a knot}
\lbl{sub.coloredj}

The {\em colored Jones function} 
$$J_K\co\BN \longto \BZ[q^{\pm}]$$
of a knot $K$ in 3--space is a sequence of Laurent polynomials,
whose $n$th term $J_K(n)$ is the Jones polynomial of a knot colored with the 
$n$--dimensional irreducible representation of $\mathfrak{sl}_2$;
see \cite{Tu}. We will normalize it by $J_{\text{unknot}}(n)=1$ for all $n$,
and (for those who worry about framings), we will assume that $K$ is 
zero-framed. 

The first two terms of the colored Jones function of a knot $K$ are 
better known. Indeed, $J_K(1)=1$, and  $J_K(2)$ coincides with the 
{\em Jones polynomial} of a knot $K$, defined by Jones in \cite{J}.
Although we will not use it, note that the colored Jones function of a knot 
essentially encodes the Jones polynomial of a knot and its connected parallels.

The starting point for our paper is the key property that 
the colored Jones function is $q$--holonomic, as was shown 
in joint work with TTQ Le; see \cite{GL}. 
Informally, a $q$--holonomic function is one that satisfies a nontrivial
linear recursion relation, with appropriate coefficients. A convenient way 
to describe recursion relations is the {\em operator point of view} which we 
now describe.

\subsection{The characteristic variety of a knot}
\lbl{sub.characteristic}

Consider the ring $\F$ of {\em discrete functions} $f\co \BN \longto \BQ(q)$,
and define the linear operators $E$ and $Q$ on $\F$ which act on a discrete 
function $f$ by:
$$(Q f)(n)=q^{n}f(n) \hspace{1cm} (E f)(n) = f(n+1).$$
It is easy to see that $EQ=qQE$, and that $E,Q$ generate a noncommutative 
{\em Weyl algebra} (often called a $q$--Weyl algebra) with presentation
$$\A=\BZ[q^{\pm}]\la Q,E \ra/(EQ=qQE).$$ 

Given a discrete function $f$, consider the set 
$$\calI_f=\{P \in \A \, | Pf=0\}.$$

It is easy to see that $\calI_f$ is a left ideal of the Weyl algebra,
the so-called {\em recursion ideal} of $f$.

If $P \in \calI_f$, we may think of the equation $P f=0$ as a {\em linear 
recursion relation} on $f$. Thus, the set of linear recursion relations that 
$f$ satisfies may be identified with the recursion ideal $\calI_f$.

\begin{definition}
\lbl{def.qholo} 
We say that $f$ is $q$--{\em holonomic} iff $\calI_f \neq 0$. In other words,
a discrete function is $q$--holonomic iff it satisfies a nontrivial linear 
recursion relation. 
\end{definition}

Consider the quotient $\B=\BZ[E,Q]$ of the Weyl algebra and let 
\begin{equation}
\lbl{eq.epsilon}
\e\co\A\longto\B
\end{equation} 
be the {\em evaluation map} at $q=1$.

\begin{definition}
\lbl{def.chara}
If $I$ is a left ideal in $\A$, we define its
{\em characteristic variety} $\ch(I) \sub (\BC^{\star})^2$ by
$$\ch(I)=\{(x,y) \in (\BC^{\star})^2 | P(x,y)=0 \,\,\, \text{for all} \,\,\,
  P \in \e(I) \}.$$
If $f$ is a $q$--holonomic function, then we define its {\em characteristic 
variety} to be $\ch(\calI_f)$. 
Finally, if $K$ is a knot in 3--space, we define its {\em characteristic
variety} $\ch(K)$ to be $\ch(J_K)$.
\end{definition}

We will make little distinction between a variety $V \sub (\BC^{\star})^2$ and
its closure $\overline{V} \sub \BC^2$. For those proficient in holonomic 
functions, please note that our definition of characteristic variety
does {\em not} agree with the one
commonly used in holonomic functions. The latter uses only the symbol (ie
the leading $E$--term) of recursion relations.

\subsection{The deformation variety of a knot}
\lbl{sub.deformation}

The deformation variety of a knot is the character variety of
$\SL_2(\BC)$ representations of the knot complement, viewed from their
restriction to the boundary torus. The deformation variety of a knot is
of fundamental importance to hyperbolic geometry, and to geometrization,
and was studied extensively by Cooper et al and Thurston; see \cite{CCGLS}
and \cite{Th}. 

Given a knot $K$ in $S^3$, consider the complement $M=S^3-\text{nbd}(K)$
(a 3--manifold with torus boundary $\pt M \cong T^2$), and 
the set 
$$R(M)=\Hom(\pi_1(M),\SL_2(\BC))$$ 
of representations of $\pi_1(M)$ into $\SL_2(\BC)$. This has the structure
of an affine algebraic variety defined over $\BQ$, 
on which $\SL_2(\BC)$ acts by conjugation
on representations. Let $X(M)$ denote the algebrogeometric quotient. There
is a natural restriction map $X(M) \longto X(\pt M)$, induced by the
inclusion $\pt M \sub M$. Notice that 
$\pi_1(\pt M)\cong \BZ^2$, generated by a meridian and longitude of $K$.
Restricting attention to representations of $\pi_1(\pt M)$ which are
upper diagonal, we may identify the character variety of $\pt M$
with $(\BC^2)^\star$,
parametrized by $L$ and $M$, the upper left entry of meridian and longitude.
Cooper, Culler, Gillet, Long and Shalen \cite{CCGLS} define the
{\em deformation variety} $D(K)$ to be the image of $X(\pt M)$ in
$(\BC^\star)^2$.

\subsection{The conjecture}
\lbl{sub.conj}

Recall that every affine subvariety $V$ in $\BC^2$ is the disjoint union
$V_0 \sqcup V_1 \sqcup V_2$ where $V_i$ is a subvariety of $V$ of {\em pure 
dimension} $i$. 

We say that two algebraic subvarieties $V$ and $V'$ of $\BC^2$ are 
{\em essentially equal} iff $V_1$ is equal to $V'_1$ union some $y$--lines,
where a $y$--line in $\BC^2=\{(x,y) \, | x,y \in \BC\}$ is a line $y=a$ for 
some $a$.

\begin{conjecture}[The Characteristic equals Deformation Variety Conjecture]
\lbl{conj.1}
For every knot in $S^3$, the characteristic and  deformation varieties 
are essentially equal.
\end{conjecture}

Questions similar to the above conjecture and its polynomial version
(Conjecture \ref{conj.2} below) were also raised by Frohman and Gelca 
who studied the colored Jones function of a knot via {\em Kauffman bracket
skein theory}, \cite{Ge}. Our approach to recursion relations in \cite{GL}
and here is via statistical mechanics sums and holonomic functions.

A modest corollary of the above conjecture is the following:

\begin{corollary}
\lbl{cor.1}
If a knot has nontrivial deformation variety (eg. the knot is 
hyperbolic), then it has nontrivial colored Jones function.
\end{corollary}

\begin{remark}
\lbl{rem.Atrivial}
Despite our improved understanding of the geometry of 3--mani\-folds,
it is unknown at present whether the deformation variety of a knot
complement is positive dimensional. If a knot is hyperbolic or torus,
then it is, by above mentioned work of Thurston and Cooper et al.
If a knot is a satellite, then it is not known, due to the presence
of {\em forbidden representations}, explained by Cooper--Long in 
\cite[section 9]{CL}.
\end{remark}

As evidence for the conjecture, we will show by a direct calculation that:

\begin{proposition}
\lbl{prop.1}
Conjectures \ref{conj.1} and \ref{conj.2} are true for the trefoil and 
figure--8 knots.
\end{proposition}

Let us end this section with three comments:

\begin{remark}
\lbl{rem.1}
Conjecture \ref{conj.1} may be translated as an equality of two polynomials
with two commuting variables and integer coefficients; see Conjecture 
\ref{conj.2} below. Since these polynomials are computable by elimination,
it follows that Conjecture \ref{conj.1} is in principle a decidable question.
This is in contrast to the {\em Hyperbolic Volume Conjecture} (due to 
Kashaev--Murakami--Murakami; see \cite{Ka,MM}) which involves
the existence and identification of a limit of complex numbers.
\end{remark}

\begin{remark}
\lbl{rem.volume}
Both Conjecture \ref{conj.1} and the Hyperbolic Volume Conjecture state
a relationship between the colored Jones function of a knot and hyperbolic
geometry. 
Combining both conjectures, it follows that the colored Jones function of
a hyperbolic knot determines the volume of the hyperbolic 3--manifolds obtained
by Dehn surgery on the knot. Indeed, the variation of the volume function
depends on the restriction of a path of $\SL_2(\BC)$ representations to
the boundary of the knot complement. Furthermore, the polynomial that
defines the deformation variety
can compute the variation of the volume function; see Cooper et al
\cite[section 4.5]{CCGLS} and also Yoshida \cite{Y} and 
Neumann--Zagier \cite[equation (47)]{NZ}.
\end{remark}

\begin{remark}
\lbl{rem.2}
Conjecture \ref{conj.1} reveals a close relation between the colored Jones
function of a knot and its deformation variety. It does not explain though
why we ought to look at characters of $\SL_2(\BC)$ representations. There
is a generalization to higher rank groups, which we present in Section 
\ref{sec.high}. We warn the reader that there is no evidence for this 
generalization.
\end{remark}

\subsection{Acknowledgement}
An earlier version of this paper circulated informally at the CassonFest
in the University of Texas at Austin in May 2003. The author wishes to thank
the organizers, C Gordon, J Luecke and A Reid for arranging a stimulating
conference and for their hospitality.

In addition, the author wishes to thank A Riese and D Thurston for 
enlightening conversations, and the anonymous referee.

The author was supported in part by NSF and BSF.

\section{A polynomial version of Conjecture \ref{conj.1}}
\lbl{sec.reformulation}

\subsection{The $A$--polynomial of a knot}
\lbl{sub.A}

Recall the definition of the deformation variery of a knot from Section 
\ref{sub.deformation}. 
Since projection of affine algebraic varieties corresponds to elimination
in their corresponding ideals (see \cite{CLO}), it is clear that the
deformation variety of a knot can in principle be computed via elimination.

In fact, according to \cite{CCGLS}, the deformation variety $D(K)$ of a knot
$K$ is essentially equal to a complex curve in $\BC^2$ which is defined
by the zero-locus of the so-called $A$--{\em polynomial} $A(K)$ of $K$, where
the latter lies in $\BZ[L,M^2]$.
Here $A$ stands for {\em affine} and not for Alexander.

\subsection{A noncommutative version of the $A$--polynomial}
\lbl{sub.Aq}

In this section we 
define a noncommutative version of the $A$--polynomial of a knot.

If the Weyl algebra $\A$ were a principal ideal domain, every left ideal
(such as the recursion ideal of a discrete function) would be generated by a 
polynomial in noncommuting variables $E$ and $Q$. This polynomial would
be the noncommutative $A$--polynomial of an ideal. Applying this to the
recursion ideal of $J_K$
would allow us to define the noncommutative $A$--polynomial of a knot.

Unfortunately, the algebra $\A$ is not a principal ideal domain. One way to 
get around this problem is to invert polynomials in $Q$, as we now explain.
Consider the  {\em Ore algebra}
$\A_{\loc}=\BK[E,\s]$ over the field $\BK=\BQ(q,Q)$, where
$\s$ is the automorphism of $\BK$ given by 
\begin{equation}
\lbl{eq.sigma}
\s(f)(q,Q)=f(q,qQ).
\end{equation}
Additively, we have
$$\A_{\loc} = \left\{\sum_{k=0}^\infty a_k E^k \qua \left|\qua a_k \in \BK,
  \qua a_k=0 ~\text{for}~ k \gg 0 \right. \right\},$$ 
where the multiplication of monomials given by $a E^k \cdot b E^l=a \s^k(b) 
E^{k+l}$.

Recall the ring $\F$ of discrete functions $f\co \BN \longto \BQ(q)$, and its
quotient ring $\tilde{\F}$ under the equivalence relation $f \sim g$ iff
$f(n)=g(n)$ for all but finitely many $n$. Then, $\A_{\loc}$ acts on 
$\tilde{\F}$. In particular, if $f$ is a discrete function, we may define
its recursion ideal, with respect to $\A_{\loc}$. We will call $f$ 
$q$--holonomic with respect to $\A_{\loc}$ iff its recursion ideal with
respect to $\A_{\loc}$ does not vanish.

By clearing out denominators, it is easy to see that if $f$ is a discrete
function, then it is $q$--holonomic with respect to $\A$ iff it is
$q$--holonomic with respect to $\A_{\loc}$.

It turns out that
every left ideal in $\A_{\loc}$ is {\em principal}; see \cite[chapter 2, 
exercise 4.5]{Cou}. Given a left ideal $I$ of $\A_{\loc}$, let $A_q(I)$ 
denote a generator of $I$, with the following properties:
\begin{itemize}
\item
$A_q(I)$ has smallest $E$--degree and lies in $\A$.
\item
We can write $A_q(I)=\sum_k a_k E^k$ where $a_k \in \BZ[q,Q]$
are coprime (this makes sense since $\BZ[q,Q]$ is a unique factorization 
domain). 
\end{itemize}
These properties uniquely determine $A_q(I)$ up to left multiplication
by $\pm q^a Q^b$ for integers $a,b$.

\begin{definition}
\lbl{def.principal}
Given a left ideal $I$ in $\A$, we define its $A_q$--{\em polynomial} 
$A_q(I) \in \A$ to be $A_q(I)$. 
Given a knot $K$ in $S^3$, we define its 
$A_q$--polynomial $A_q(K)$ to be the $A_q$--polynomial of the 
$\A_{\loc}$--recursion ideal of $J_K$.
\end{definition}

Recall from Section \ref{sub.A} that the $A$ polynomial of a knot lies
in the ring $\BZ[L,M^2]$ which we will identify with
$\BZ[E,Q]$ by $L=E$ and $M=Q^{1/2}$. In other words, 

\begin{definition}
\lbl{def.identify}
We identify the geometric pair $(L,M^2)$ of ({\em meridian, longitude}) 
of a knot $K$ with the pair $(E,Q)$ of {\em basic operators} which act
on the colored Jones function of $K$.
\end{definition}

Let us comment on this definition. It is not too surprising
that the meridian variable $M$ is identified with $Q$, the multiplication
by $q^n$. This is foreshadowed by the {\em Euler expansion} of the colored 
Jones function in terms of powers of $q^n$ and $q-1$, \cite{G}. The physical
meaning of this expansion is, according to Rozansky, a  Feynman diagram 
expansion around a $U(1)$--connection in the knot complement with holonomy 
$q^n$, \cite{R}. Thus, it is not surprising that $M^2=Q$.

It is more surprising that the longitude variable $L$ corresponds 
to the shift operator $E$. This can be explained in the following way. 
According to Witten (see \cite{Wi}), the Jones polynomial $J_K(n)$ 
of a knot $K$ is the average over an infinite dimensional space of connections,
of the {\em trace of the holonomy around $K$}, where the trace is computed
in the $n$--dimensional representation of $\mathfrak{sl}_2$. To a leading
order term, computing traces in the $n$--dimensional representation is
equivalent to computing traces of an $(n-1,1)$--connected parallel of the knot
in the $2$--dimensional representation. Thus, increasing $n$ by $1$ corresponds
to going once more around the knot. Since holonomy and longitude are 
synonymous notions, this explains in some sense the relation $E=L$.

\begin{conjecture}
\lbl{conj.2}{\rm(The AJ Conjecture)}\footnote{AJ are the initials of the 
$A$--polynomial and the colored Jones polynomial}\qua
For every knot in $S^3$, $A(K)(L,M)=\e A_q(K)(L,M^2)$.
\end{conjecture}

\begin{lemma}
\lbl{lem.conj}
Conjecture \ref{conj.2} implies Conjecture \ref{conj.1}.
\end{lemma}

\begin{proof}
Consider $f,g \in \BZ[E,Q]$. Let us say that $f$ is {\em essentially equal} to
$q$ if their images in $\BQ(Q)[E]$ are equal.
In other words, $f$ is essentially equal to $g$ iff $f/g$ is a rational
function of $Q$.

If $V(f)=\{(x,y) \in \BC^2 \, | f(x,y)=0 \}$ denotes the variety of zeros
of $f$, then it is easy to see that if $f$ is essentially equal to $g$, then
$V(f)$ is essentially equal to $V(g)$.

It is easy to see that the characteristic (resp. deformation) variety is 
essentially equal to $V(\e A_q)$ (resp. $V(A)$). The result follows.
\end{proof}

\begin{remark}
\lbl{rem.symmetry}
Conjecture \ref{conj.2} is consistent with the behavior of the colored Jones 
function and the $A$--polynomial under mirror image, changing the orientation 
of the knot, and $\BZ_2$--symmetry. For the behavior of the $A$--polynomial
under these operations, see Cooper--Long: \cite[proposition 4.2]{CL}.
On the other hand, the colored Jones function satisfies the symmetry
$J(n)=J(-n)$. Moreover, $J$ is invariant under the change of orientation of a
knot and changes under $q\to q^{-1}$ under mirror image.
\end{remark}

\subsection{Computing the $A_q$ polynomial of a knot}
\lbl{sub.computeAq}

Section \ref{sec.reformulation} defines the $A_q$ polynomial of a knot $K$.
This section explains how to compute the $A_q$ polynomial of a knot. For more 
details, we refer the reader to \cite{GL}.

Starting from a generic planar projection of a knot $K$, 
it was shown in \cite[section 3.2]{GL} that the colored Jones
function of a knot $K$ can be written as a {\em multisum}
\begin{equation}
\lbl{eq.multi}
J_K(n)=\sum_{k_1,\dots,k_r=0}^\infty F(n,k_1,\dots,k_r)
\end{equation}
of a {\em proper $q$--hypergeometric function} $F(n,k_1,\dots,k_r)$.
For a fixed positive $n$, only finitely many terms are nonzero. 
Of course, $F$ depends on a planar projection of $K$. The key property
is that $F$ is $q$--holonomic in all $r+1$ variables, and that it follows from
first principles that multisums of $q$--holonomic functions are $q$--holonomic
in all remaining free variables.

Working with the Weyl algebra $\A_r$ of $r+1$ variables, and using
the fact that $F$ is $q$--proper hypergeometric, we may write
$EF/F=A/B$ and $E_iF/F=A_i/B_i$ for polynomials $A,B, A_i, B_i \in
\BQ(q)[q^n,q^{k_1},\dots,q^{k_r}]$. Replacing $q^n$ by $Q$ and $q^{k_i}$
by $Q_i$, it follows that the recursion ideal of $F$ in the Weyl
algebra $\A_{r+1}$ is generated by $BE-A, B_1 E_1-A_1,\dots, B_r E_r-A_r$.

The creative telescoping method of Wilf--Zeilberger (the so-called
{\em WZ algorithm}) produces from these 
generators of $F$, via noncommutative elimination, operators that annihilate 
$J_K$. For a discussion of Wilf--Zeilberger's algorithm, see \cite{Z,WZ,PWZ}
and also \cite[section 5]{GL}. For an implementation of the algorithm, see
\cite{PR1,PR2}.

Applying the WZ algorithm to equation \eqref{eq.multi}, we are guaranteed
to get an operator $P \in \A_{\loc}$ such that $P J_K=0$.
It follows that $A_q(K)$ is a right divisor of $P$. In other words, there
exist an operator $P_1 \in \A_{\loc}$ such that $P_1P=A_q(K)$.
We caution however that the WZ algorithm does not give in general a minimal 
order difference operator. For a thorough discussion of this matter, see
\cite[p164]{PWZ}. In other words, $P$ need not equal to $A_q(K)$.

The problem of computing right factors of an operator has been solved in
theory by Petkov\v sek in \cite{BP}. A computer implementation of this solution
is not available at present.

In case we are looking for right factors of degree $1$ (this is equivalent
to deciding whether a discrete function has closed form), there is an 
algorithm {\tt qHyper} of Petkov\v sek which decides about this problem
in real time; see \cite{PWZ}.

In the special examples that we will consider, namely the colored Jones
function of $3_1$ and $4_1$ knots, we can bypass the thorny issue of right 
factorization of an operator.

\section{Proof of the conjecture for the trefoil and figure--8 knots}
\lbl{sec.proof}

\subsection{The colored Jones function and the $A$--polynomial
of the $3_1$ and $4_1$ knots}
\lbl{sub.jj}

Habiro \cite{H} and Le give the following formula for the colored Jones 
function
of the left-handed trefoil ($3_1$) and figure--8 ($4_1$) knots:
\begin{eqnarray}
\lbl{eq.trefoil}
J_{3_1}(n) &=& \sum_{k=0}^\infty (-1)^k q^{k(k+3)/2} 
q^{nk} (q^{-n-1};q^{-1})_k (q^{-n+1};q)_k \\
\lbl{eq.figure8}
J_{4_1}(n) &=& \sum_{k=0}^\infty q^{nk} (q^{-n-1};q^{-1})_k (q^{-n+1};q)_k.
\end{eqnarray}
where we define the {\em rising} and {\em falling factorials} for $k > 0$
by:
\begin{eqnarray*}
(a;q)_k &=& (1-a)(1-aq) \dots (1-aq^{k-1}) \\
(a;q^{-1})_k &=& (1-a)(1-aq^{-1}) \dots (1-aq^{-k+1}) 
\end{eqnarray*}
and $(a;q)_0=(a;q)_0=1$. Notice that the sums in in equations 
\eqref{eq.trefoil} and \eqref{eq.figure8} have compact support, namely
for each positive $n$, only the terms with $k \leq n$ contribute.

These formulas are discussed in detail in Masbaum \cite[theorem 5.1]{Ma}, in 
relation to the cyclotomic expansion of the colored Jones function of twist 
knots. To compare Masbaum's formula with the one given above, keep in
mind that:
\begin{eqnarray*}
S(n,k)&:=& q^{nk} (q^{-n-1};q^{-1})_k (q^{-n+1};q)_k \\
&=& 
\frac{\{ n-k\}\{n-k+1\} \dots
\{n+k\}}{\{n\}} \\
&=& \prod_{j=1}^k ((q^{n/2}-q^{-n/2})^2-(q^{j/2}-q^{-j/2})^2)
\end{eqnarray*}
where $\{m\}=q^{m/2}-q^{-m/2}$.

On the other hand, Cooper et al \cite{CCGLS} compute the $A$--polynomial of the
$3_1$ and $4_1$ knots, as follows:
\begin{eqnarray}
\lbl{eq.Atrefoil}
A(3_1) &=&(L-1)( L+M^6) \\
\lbl{eq.Afigure8}
A(4_1) &=&(L-1)( -L+LM^2+M^4+2LM^4+L^2M^4+LM^6-LM^8)
\end{eqnarray}
where we include the factor $L-1$ in the $A$--polynomial which corresponds
to the abelian representations of the knot complement.

\subsection{Computer calculations}
\lbl{sec.computer}

The colored Jones function of the $3_1$ and $4_1$ knots given in Equations
\eqref{eq.trefoil} and \eqref{eq.figure8} has no {\em closed form}. However,
it is {\em guaranteed} to obey nontrivial recursion relations. Moreover, these
relations can be found by computer.
There are various programs that can compute the recursion relations
for multisums. In maple, one may use {\tt qEKHAD} developed by
Zeilberger \cite{PWZ}. In Mathematica, one may use {\tt qZeil.m} developed by
Paule and Riese \cite{PR1,PR2}. 
We will give explicit examples in Mathematica, using Paule and Riese's {\tt
qZeil.m} package.

We start in computer talk by loading the packages:

{\footnotesize\parskip0pt
\begin{verbatim}
Mathematica 5.0 for Sun Solaris
Copyright 1988-2000 Wolfram Research, Inc.
 -- Motif graphics initialized --
In[1]:=<< qZeil.m
\end{verbatim}
\verb+q-Zeilberger Package by Axel Riese --+ \copyright \verb+RISC Linz -- V 2.35 (04/29/03)+
}
{\footnotesize\baselineskip0pt
\begin{verbatim}
In[2]:= << qMultiSum.m
\end{verbatim}
\verb+qMultiSum Package by Axel Riese --+ \copyright \verb+RISC Linz -- V 2.45 (04/02/03)+
}

Let us type the colored Jones function $J_{3_1}$ from Equation 
\eqref{eq.trefoil}:

{\footnotesize\baselineskip0pt
\begin{verbatim}
In[3]:= summandtrefoil = (-1)^k q^(k(k + 3)/2) q^(n  k) qfac[q^(-n - 1),
      q^(-1), k] qfac[q^(-n + 1), q, k]

            k  (k (3 + k))/2 + k n              -1 - n  1
Out[3]= (-1)  q                    qPochhammer[q      , -, k]
                                                        q

                  1 - n
>    qPochhammer[q     , q, k]
\end{verbatim}
}

We now ask for a recursion relation for $J_{3_1}$:

{\footnotesize\baselineskip0pt
\begin{verbatim}
In[4]:= qZeil[summandtrefoil, {k, 0, Infinity}, n, 1]

qZeil::natbounds: Assuming appropriate convergence.

                   -2 + n        2 n     -1 + 3 n       -1 + n
                  q       (-q + q   )   q         (1 - q      ) SUM[-1 + n]
Out[4]= SUM[n] == ------------------- - -----------------------------------
                              n                            n
                        -1 + q                        1 - q
\end{verbatim}
}

In other words, for $J(n)=J_{3_1}(n)$ we have: 
$$
J(n)=q^{-2+n}\frac{-q+q^{2n}}{-1+q^n} -q^{-1+3n}\frac{1-q^{-1+n}}{1-q^n} 
J(n-1),
$$
The above relation is a first order inhomogeneous recursion relation.
We may convert it into a second order homogeneous recursion relation as 
follows:

{\footnotesize\baselineskip0pt
\begin{verbatim}
In[5]:= rec31 = MakeHomRec[%, SUM[n]]

         -1 + 2 n   2    n
        q         (q  - q ) SUM[-2 + n]
Out[5]= ------------------------------- +
                    3    2 n
                   q  - q

             n        n    4    4 n    3 + n    2 + 2 n    3 + 2 n
>     ((q - q ) (q + q ) (q  + q    - q      + q        - q        -

            1 + 3 n                   n       2 n    3    2 n
>          q       ) SUM[-1 + n]) / (q  (q - q   ) (q  - q   )) +

       2 - n        n
      q      (-1 + q ) SUM[n]
>     ----------------------- == 0
                  2 n
             q - q
\end{verbatim}
}

Perhaps the reader is displeased to see the above recursion relation written
in {\em backwards shifts}, ie, {\tt SUM[-k+n]} where $k \geq 0$.
This can be converted into a recursion relation using {\em forward shifts}
by:

{\footnotesize\baselineskip0pt
\begin{verbatim}
In[6]:= ForwardShifts[% ]

         3 + 2 n   2    2 + n
        q        (q  - q     ) SUM[n]
Out[6]= ----------------------------- +
                 3    4 + 2 n
                q  - q

        -2 - n       2 + n        2 + n
>     (q       (q - q     ) (q + q     )

           4    5 + n    6 + 2 n    7 + 2 n    7 + 3 n    8 + 4 n
>        (q  - q      + q        - q        - q        + q       ) SUM[1 + n])

                                                     2 + n
                  4 + 2 n    3    4 + 2 n     (-1 + q     ) SUM[2 + n]
>        / ((q - q       ) (q  - q       )) + ------------------------ == 0
                                                  n       4 + 2 n
                                                 q  (q - q       )
\end{verbatim}
}

The next command converts the recursion relation {\tt rec31} into an
operator, where (due to {\tt Mathematica} annoyance), we use the symbol $X$ to
denote the shift $E$:

{\footnotesize\baselineskip0pt
\begin{verbatim}
In[7]:= ToqHyper[rec31[[1]] - rec31[[2]]] /. {SUM[N] -> 1, SUM[N q^c_.] :> X^c}
  /. N -> Q

         2               2       2
        q  (-1 + Q)    (q  - Q) Q
Out[7]= ----------- + -------------- +
                2         3    2   2
        Q (q - Q )    q (q  - Q ) X

                       4    3      2  2    3  2      3    4
     (q - Q) (q + Q) (q  - q  Q + q  Q  - q  Q  - q Q  + Q )
>    -------------------------------------------------------
                             2    3    2
                     Q (q - Q ) (q  - Q ) X
\end{verbatim}
}

This operator right-divides the $A_q$ polynomial of the $3_1$ knot.
Let us assume for now that it equals to the $A_q$ polynomial, after clearing
denominators.
Setting $q=1$, and replacing $X$ by $L$ and $Q$ by $M^2$, and obtain:

{\footnotesize\baselineskip0pt
\begin{verbatim}
In[8]:= Factor[ToqHyper[rec31[[1]] - rec31[[2]]] /. {SUM[N] -> 1,
          SUM[N q^c_.] :> X^c} /. {N -> Q, q -> 1}] /. {Q -> M^2, X -> L}

                         6
          (-1 + L) (L + M )
Out[8]= -(-----------------)
            2  2       2
           L  M  (1 + M )
\end{verbatim}
}

The result agrees, up to multiplication by a rational function of $M$ and a
power of $E$, with the $A$--polynomial of $3_1$ from \eqref{eq.Atrefoil}.

It remains to prove that $\mathrm{rec31}$:={\tt Out[7]} coincides with 
$A_q(3_1)$, after
clearing denominators. Notice that $\mathrm{rec31}=PA_q(3_1)$ for some operator
$P$ and $\ord_E(\mathrm{rec31})=2$, where $\ord_E(P)$ denotes the $E$--order 
of an operator $E$. Thus $\ord_E(A_q(3_1))$ is $1$ or $2$.
If $\ord_E(A_q(3_1))=1$, then $J_{3_1}$ would have a closed form. This 
problem can be decided by computer using {\tt qHyper} (see \cite{PWZ}), which 
indeed confirms that $J_{3_1}$ does not have closed form. 
Thus $\ord_E(A_q(3_1))=2=\ord_E(\mathrm{rec31})$. It follows that 
(up to left multiplication by units), $A_q(3_1)$ equals to $\mathrm{rec31}$. 
This completes the proof in the case of the trefoil.

Now, let us repeat the process for the colored Jones function of the figure 8
knot, given in Equation \eqref{eq.figure8}.

{\footnotesize\baselineskip0pt
\begin{verbatim}
In[9]:= summandfigure8 = q^(n  k) qfac[q^(-n - 1), q^(-1), k] qfac[q^(-n + 1),
  q, k]

         k n              -1 - n  1                  1 - n
Out[9]= q    qPochhammer[q      , -, k] qPochhammer[q     , q, k]
                                  q

In[10]:= qZeil[summandfigure8, {k, 0, Infinity}, n, 2]

qZeil::natbounds: Assuming appropriate convergence.

                    -1 - n       n         2 n
                   q       (q + q ) (-q + q   )
Out[10]= SUM[n] == ---------------------------- -
                                   n
                             -1 + q

            -2 + n        -1 + 2 n
      (1 - q      ) (1 - q        ) SUM[-2 + n]
>     ----------------------------------------- +
                    n        -3 + 2 n
              (1 - q ) (1 - q        )

        -2 - 2 n       -1 + n 2       -1 + n
>     (q         (1 - q      )  (1 + q      )

           4    4 n    3 + n    1 + 2 n    3 + 2 n    1 + 3 n
>        (q  + q    - q      - q        - q        - q       ) SUM[-1 + n]) /

              n        -3 + 2 n
>      ((1 - q ) (1 - q        ))
\end{verbatim}
}

gives a second-order inhomogeneous recursion relation, which we convert
into a third-order homogeneous recursion relation:

{\footnotesize\baselineskip0pt
\begin{verbatim}
In[11]:= rec41 = MakeHomRec[%, SUM[n]]

          2 + n    3    n
         q      (-q  + q ) SUM[-3 + n]
Out[11]= ----------------------------- -
              2    n     5    2 n
            (q  + q ) (-q  + q   )

        -2 - n   2    n    8    4 n      6 + n    7 + n    3 + 2 n
>     (q       (q  - q ) (q  + q    - 2 q      + q      - q        +

            4 + 2 n    5 + 2 n    1 + 3 n      2 + 3 n
>          q        - q        + q        - 2 q       ) SUM[-2 + n]) /

              n    5    2 n
>      ((q + q ) (q  - q   )) +

        -1 - n        n    4    4 n    2 + n      3 + n    1 + 2 n
>     (q       (-q + q ) (q  + q    + q      - 2 q      - q        +

            2 + 2 n    3 + 2 n      1 + 3 n    2 + 3 n
>          q        - q        - 2 q        + q       ) SUM[-1 + n]) /

                                  1 + n        n
          2    n         2 n     q      (-1 + q ) SUM[n]
>      ((q  + q ) (-q + q   )) + ----------------------- == 0
                                         n        2 n
                                   (q + q ) (q - q   )
\end{verbatim}
}

In forward shifts, we have:

{\footnotesize\baselineskip0pt
\begin{verbatim}
In[12]:= ForwardShifts[%]

           5 + n    3    3 + n
          q      (-q  + q     ) SUM[n]
Out[12]= ------------------------------ -
           2    3 + n     5    6 + 2 n
         (q  + q     ) (-q  + q       )

        -5 - n   2    3 + n    8      9 + n    10 + n    9 + 2 n
>     (q       (q  - q     ) (q  - 2 q      + q       - q        +

            10 + 2 n    11 + 2 n    10 + 3 n      11 + 3 n    12 + 4 n
>          q         - q         + q         - 2 q         + q        )

                              3 + n    5    6 + 2 n
>        SUM[1 + n]) / ((q + q     ) (q  - q       )) +

        -4 - n        3 + n    4    5 + n      6 + n    7 + 2 n    8 + 2 n
>     (q       (-q + q     ) (q  + q      - 2 q      - q        + q        -

            9 + 2 n      10 + 3 n    11 + 3 n    12 + 4 n
>          q        - 2 q         + q         + q        ) SUM[2 + n]) /

                                          4 + n        3 + n
          2    3 + n         6 + 2 n     q      (-1 + q     ) SUM[3 + n]
>      ((q  + q     ) (-q + q       )) + ------------------------------- == 0
                                                 3 + n        6 + 2 n
                                           (q + q     ) (q - q       )
\end{verbatim}
}

In operator form, {\tt rec41} becomes:

{\footnotesize\baselineskip0pt
\begin{verbatim}
In[13]:= ToqHyper[rec41[[1]] - rec41[[2]]] /. {SUM[N] -> 1, SUM[N q^c_.] :>
  X^c} /. N -> Q

                                 2      3
           q (-1 + Q) Q         q  Q (-q  + Q)
Out[13]= ---------------- + ---------------------- -
                       2      2         5    2   3
         (q + Q) (q - Q )   (q  + Q) (-q  + Q ) X

        2        8      6      7      3  2    4  2    5  2      3      2  3
>    ((q  - Q) (q  - 2 q  Q + q  Q - q  Q  + q  Q  - q  Q  + q Q  - 2 q  Q  +

           4       2             5    2   2
>         Q )) / (q  Q (q + Q) (q  - Q ) X ) +

                 4    2        3        2    2  2    3  2        3    2  3
>    ((-q + Q) (q  + q  Q - 2 q  Q - q Q  + q  Q  - q  Q  - 2 q Q  + q  Q  +

           4            2             2
>         Q )) / (q Q (q  + Q) (-q + Q ) X)
\end{verbatim}
}

where $X=E$. Let us assume that this coincides with $A_q(4_1)$, after
we clear denominators. 
Setting $q=1$, and replacing $X$ by $L$ and $Q$ by $M^2$, and 
obtain:

{\footnotesize\baselineskip0pt
\begin{verbatim}
In[14]:= Factor[ToqHyper[rec41[[1]] - rec41[[2]]] /. {SUM[N] -> 1,
          SUM[N q^c_.] :> X^c} /. {N -> Q , q -> 1}] /. {Q -> M^2, X -> L}

                          2    4        4    2  4      6      8
         (-1 + L) (L - L M  - M  - 2 L M  - L  M  - L M  + L M )
Out[14]= -------------------------------------------------------
                              3  2       2 2
                             L  M  (1 + M )
\end{verbatim}
}

The result agrees, up to multiplication by a rational function of $M$ and a
power of $E$, with the $A$--polynomial of $4_1$ from \eqref{eq.Afigure8}.

It remains to prove that $\mathrm{rec41}$:={\tt Out[13]} is equal,
up to units, to $A_q(4_1)$. 
Notice that $\mathrm{rec41}=PA_q(4_1)$ for some operator
$P$ and $\ord_E(\mathrm{rec41})=3$. 
Thus $\ord_E(A_q(4_1))$ is $1$ or $2$ or $3$.

If $\ord_E(A_q(4_1))=1$, then $J_{4_1}$ would have a closed form. This 
problem can be decided by computer using {\tt qHyper} (see \cite{PWZ}), which 
indeed confirms that $J_{3_1}$ does not have closed form. 

If $\ord_E(A_q(4_1))=2$, recall the map $\e$ which evaluates at $q=1$.
We have: $\e\mathrm{rec41}=\e P \, \e A_q(4_1)$. Since 
$\ord_E(\e\mathrm{rec41})=3$, it follows that we must have 
$\ord_E(\e A_q(4_1))=2$.

Furthermore, the computer calculation above shows that $\e A_q(4_1)$
divides $A(4_1)$. The latter, given by Equation \eqref{eq.Afigure8} can
be factored as a product of two irreducible polynomials of $E$--degree
$1$ and $2$.

On the other hand, Lemma \ref{lem.divideQ} below implies that $E-1$ divides
$(\e A_q(4_1))|_{Q=1}$. Combining these facts, it follows that $\e A_q(4_1)=
A(4_1)$ (and therefore, also $A_q(4_1)$) is of $E$--degree $3$, a contradiction 
to our hypothesis.

Thus, it follows that $\ord_E(A_q(4_1))=3=\ord_E(\mathrm{rec41})$. This
implies that, up to left multiplication by units, $A_q(4_1)$ coincides with
$\mathrm{rec41}$. This concludes the proof in the case of the figure--8 knot.


\begin{lemma}
\lbl{lem.divideQ}
For every knot $K$, $\e\A_q(K)(1,1)=0$.
\end{lemma}

\begin{proof}
Recall that the colored Jones function of a knot $K$ is given by a multisum 
formula of a $q$--proper hypergeometric function. Consider the evaluation of 
the colored Jones function $\e J_K$ at $q=1$. This is a discrete function 
which is given by a multisum of a proper hypergeometric function. Applying
the WZ algorithm, it follows that $\e_Q\e\A_q(K)$ annihilates $\e J_K$,
where $\e_Q$ is the evaluation at $Q=1$. However, $\e J_K(n)=1$ for all $n$; 
see \cite{GL}. Thus $E-1$ divides $\e_Q\e\A_q(K)$. The result follows.
\end{proof}

\section{Higher rank groups}
\lbl{sec.high}

The purpose of this section is to formulate a generalization of the
characteristic and deformation varieties of a knot to higher rank groups. 

Consider a {\em simple} simply connected compact Lie group $G$ with Lie 
algebra $\fg$ and complexified group $G_{\BC}$. Let $\La \cong \BZ^r$ 
denote its weight
lattice, which is a free abelian group of rank $r$, the rank of $G$,
and let $\La_+ \cong \BN^r$ denote the cone of positive dominant weights.

One can define the $\fg$--{\em colored Jones function} 
$$J_{\fg}\co \BN^r \longto \BZ[q^{\pm}].$$ 
In \cite{GL}, we showed that $J_{\fg}$ is $q$--holonomic, with respect to the 
Weyl algebra of $r$ variables:
$$
\A_r= \frac{\BZ[q^{\pm}]\la Q_1, \dots, Q_r, E_1, 
\dots, E_r \ra}{(\text{Rel}_q)}
$$
where the relations are given by:
\begin{equation*}
\tag{$\text{Rel}_{q}$}
\begin{aligned}
Q_i Q_j&= Q_j Q_i  & E_i E_j &= E_j E_i \\
Q_i E_j&= E_j Q_i ~ \text{for} ~ i \neq j & E_i Q_i &= q Q_i E_i
\end{aligned}
\end{equation*}

Loosely speaking, holonomicity of a discrete function of $r$ variables
means that it satisfies $r$ independent linear recursion relations.

A precise definition in several equivalent forms was given in 
\cite[section 2]{GL}. For the benefit of the reader, we recall here
the definition in its form most useful for our purposes.

Given a discrete function $f\co \BN^r\longto\BQ(q)$, we define the {\em 
recursion ideal} $\calI_f$ and the $q$--{\em Weyl module} $M_f$ by:
$$\calI_f=\{P \in \A_r \, | Pf=0\} \qquad
  M_f:=\A_r f \cong \A_r/\calI_f.$$

$M_f$ is a cyclic left $\A_r$ module. Every finitely generated left 
$\A_r$ module has a Hilbert dimension. In case $M=\A_r/I$ is
cyclic, its 
{\em Hilbert dimension} $d(M)$ is defined as follows. Let $F_m$ be
the sub-space of $\A_r$ spanned by polynomials in $Q_i,E_i$ of
total degree $\leq m$. Then the module $\A_r/I$ can be
approximated by the sequence $F_m/(F_m\cap I), m=1,2,...$. It
turns out that, for $m \gg 1$,  the dimension of the vector space 
$F_m/(F_m\cap I) \otimes_{\BZ[q^{\pm}]} \BQ(q)$ (over the field $\BQ(q)$) is a 
polynomial in $m$ of degree equal (by definition) to $d(M)$.

Bernstein's {\em famous inequality} (proved by Sabbah in the
$q$--case, \cite{Sa}) states that $d(M) \geq r$, if $M\neq 0$ and $M$ has
{\em no monomial torsions}, ie, any non-trivial element of $M$ cannot be
annihilated by a monomial in $Q_i,E_i$. Note that the  left
$\A_r$ module $M_f:=\A_r \cdot f \cong \A_r/\calI_f$ does not
have monomial torsion.

\begin{definition}
\lbl{def.qholor}  
We say that a discrete function $f$ is
$q$--holonomic if $d(M_f)\le r$.
\end{definition}

Note that if $d(M_f)\le r$, then by Bernstein's inequality, either
$M_f=0$ or $d(M_f)=r$. The former can happen only if $f=0$. Of course,
for $r=1$, definitions \ref{def.qholo} and \ref{def.qholor} agree.

Let us now define the characteristic variety of a cyclic $\A_r$ module 
$M=\A_r/I$. Let 
$$\B_r=\BZ[Q_1, \dots, Q_r,E_1, \dots, E_r]$$
and $\e\co \A_r \longto \B_r$ denote the evaluation map at $q=1$.

\begin{definition}
\lbl{def.high}
The {\em characteristic variery} $\ch(M)$ of $M$ is defined by
$$\ch(M)=
  \{(x,y) \in (\BC^\star)^{2r} ~ | ~ P(x,y)=0 ~ \text{for all} ~
  P \in \e(I \cap \A_r) \}$$
\end{definition}

This definition may be extended to define the characteristic variety of 
finitely generated left $\A_r$ modules.
As before, we will make little distinction between the characteristic variety
and its closure in $\BC^{2r}$.

\begin{lemma}
\lbl{lem.high}
If $M$ is a $q$--holonomic $\A_r$ module, then $\dim_{\BC}\ch(M)
\geq r$.
\end{lemma}

\begin{proof}
Since $M$ is $q$--holonomic, it follows that the Hilbert dimension of 
$(\A_r\otimes \BQ(q))/I$ is $r$, and from this it follows that the Hilbert
dimension of $(\A_r\otimes \BQ(q))/I$ for generic $q \in \BC$ is $r$.
Since dimension is upper semicontinuous and it coincides with the Hilbert
dimension at the generic point \cite{S}, the result follows.
\end{proof}

\begin{definition}
\lbl{def.high1}
If $K$ is a knot in $S^3$, and $G$ as above, 
we define its $G$--{\em characteristic variety} $V_G(K)
\sub \BC^{2r}$ to be the characteristic variery of its $\fg$--colored Jones
function.
\end{definition}

Similarly to the case of $\SL_2(\BC)$, given a knot $K$ in $S^3$, 
consider the complement $M=S^3-\text{nbd}(K)$ and the set $R_{G_{\BC}}(M)$ of
representations of $\pi_1(M)$ into $G_{\BC}$. This has the structure
of an affine algebraic variety, on which $G_{\BC}$ acts by conjugation
on representations. Let $X_{G_{\BC}}(M)$ denote the algebrogeometric 
quotient. There is a natural restriction map $X_{G_{\BC}}(M) \longto 
X_{G_{\BC}}(\pt M)$. Notice that 
$\pi_1(\pt M)\cong \BZ^2$, generated by the meridian and longitude of $K$.
Restricting attention to representations of $\pi_1(\pt M)$ which are
upper diagonal with respect to a Borel decomposition, we may identify the 
character variety $X_{G_{\BC}}(\pt M)$ with $T^2$ where $T$ is a maximal torus
in $G_{\BC}$. 

\begin{definition}
\lbl{def.high2}
The $G_{\BC}$--{\em deformation variety} $D_{G_{\BC}}(K)$ of $K$ is the image
of $X_{G_{\BC}}(\pt M)$ in $T^2$. 
\end{definition}

Notice that the maximal torus $T$ of $G_{\BC}$ can be identified with
$(\BC^\star)^r$, once we choose fundamental weights $\l_i$. 
This allows us to identify the values of meridian and longitude with 
$T^2$. Notice further that the deformation variety of a knot contains
an $r$--dimensional component which corresponds to abelian representations.

Let us say that two varieties $V$ and $V'$ in $\BC^{2r}=\{(x,y)| \, x,y \in
\BC^r\}$ are {\em essentially equal} if the pure $r$--dimensional part of $V$
equals to that of $V'$ union some $r$--dimensional varieties of the form
$f(y)=0$.

\begin{question}
\lbl{eq.high}
Is it true that for every $G$ as above and for every knot $K$, 
the characteristic and deformation varieties 
$V_{G}(K)$ and $D_{G_{\BC}}(K)$
are essentially equal?
\end{question}

\let\olditem\bibitem
\def\bibitem#1]#2{\olditem{#2}}

\Addresses\recd

\begin{thebibliography}{[EMSS]}

\bibitem[BP]{BP}
{\bf M Bronstein}, {\bf M Petkov\v sek},
{\em An introduction to pseudo-linear algebra},
from: ``Algorithmic complexity of algebraic and geometric models 
(Creteil, 1994)'', Theoret. Comput. Sci. 157  (1996) 3--33
\MR{1383396}

\bibitem[CCGLS]{CCGLS}
{\bf D Cooper}, {\bf D\,M Culler}, {\bf H Gillet}, {\bf D Long},
{\bf P Shalen},
{\em Plane curves associated to character varieties of $3$--manifolds}, 
Invent. Math. 118 (1994) 47--84
\MR{1288467}

\bibitem[CL]{CL}
{\bf D Cooper}, {\bf D Long},
{\em Remarks on the $A$--polynomial of a knot},
J. Knot Theory Ramifications 5 (1996) 609--628
\MR{1414090}

\bibitem[CLO]{CLO}
{\bf D Cox}, {\bf J Little}, {\bf D O'Shea},
{\em Ideals, varieties, and algorithms. An introduction to 
computational algebraic geometry and commutative algebra}
(second edition), Undergrad. texts series, Springer--Verlag (1997)
\MR{1417938}

\bibitem[Cou]{Cou}
{\bf S Coutinho},
{\em A primer of algebraic $D$--modules},
London Math. Soc. Student Texts 33, Cambridge University Press (1995)
\MR{1356713}

\bibitem[GL]{GL}
{\bf S Garoufalidis}, {\bf T\,T\,Q Le},
{\em The colored Jones function is $q$--holonomic},
\arxiv{math.GT/0309214}

\bibitem[G]{G}
{\bf S Garoufalidis}, {\bf T\,T\,Q Le},
{\em Difference and differential equations for the colored Jones function},
preprint (2003)

\bibitem[Ge]{Ge}
{\bf R Gelca},
{\em On the relation between the $A$--polynomial and the Jones polynomial}, 
Proc. Amer. Math. Soc. 130 (2002) 1235--1241
\MR{1873802}

\bibitem[H]{H}
{\bf K Habiro},
{\em On the quantum $sl_2$ invariants of knots and integral homology spheres},
\gtmref{4}{2002}{5}{55}{68}
\MR{2002603}

\bibitem[J]{J}
{\bf V Jones},
{\em Hecke algebra representation of braid groups and link polynomials}, 
Ann. of Math. 126 (1987) 335--388
\MR{0908150}

\bibitem[Ka]{Ka}
{\bf R Kashaev},
{\em The hyperbolic volume of knots from the quantum dilogarithm},
Modern Phys. Lett. A 39 (1997) 269--275
\MR{1434238}

\bibitem[Ma]{Ma}
{\bf G Masbaum},
{\em Skein-theoretical derivation of some formulas of Habiro},
\agtref3{2003}{17}{537}{556}
\MR{1997328}

\bibitem[MM]{MM}
{\bf H Murakami}, {\bf J Murakami},
{\em The colored Jones polynomials and the simplicial volume of a knot},
Acta Math.  186 (2001) 85--104
\MR{1828373}

\bibitem[NZ]{NZ}
{\bf W Neumann}, {\bf D Zagier},
{\em Volumes of hyperbolic three-manifolds},
Topology 24 (1985) 307--332
\MR{0815482}

\bibitem[PR1]{PR1}
{\bf P Paule}, {\bf A Riese},
{\em A Mathematica $q$--Analogue of Zeilberger's Algorithm Based on an 
Algebraically Motivated Approach to $q$--Hypergeometric Telescoping}, 
from: ``Special Functions, $q$--Series and Related Topics'', 
Fields Inst. Commun., 14 (1997) 179--210
\MR{1448687}

\bibitem[PR2]{PR2}
{\bf P Paule}, {\bf A Riese},
Mathematica software:\par
{\footnotesize\url{http://www.risc.uni-linz.ac.at/research/combinat/risc/software/qZeil/}}

\bibitem[PWZ]{PWZ}
{\bf M Petkov\v sek}, {\bf H\,S Wilf}, {\bf D Zeilberger},
{\em $A=B$}, 
A\,K Peters Ltd, Wellesley, MA (1996)
\MR{1379802}


\bibitem[R]{R}
{\bf L Rozansky},
{\em The universal $R$--matrix, Burau representation, and the 
Melvin--Morton expansion of the colored Jones polynomial},
Adv. Math. 134 (1998) 1--31
\MR{1612375}

\bibitem[Sa]{Sa}
{\bf C Sabbah},
{\em Syst\'emes holonomes d'\'equations aux $q$--diff\'erences},
from: ``$D$--modules and microlocal geometry, Lisbon (1990)'',
de Gruyter, Berlin (1993) 125--147
\MR{1206016}

\bibitem[S]{S}
{\bf I Shafarevich},
{\em Basic algebraic geometry I},
Springer--Verlag (1995)
\MR{1287418}

\bibitem[Th]{Th}
{\bf W Thurston},
{\em The geometry and topology of 3-manifolds},
Lecture notes, Princeton (1977)

\bibitem[Tu]{Tu}
{\bf V Turaev},
{\em The Yang--Baxter equation and invariants of links},
Invent. Math. 92 (1988) 527--553
\MR{0939474}

\bibitem[WZ]{WZ}
{\bf H Wilf}, {\bf D Zeilberger},
{\em An algorithmic proof theory for hypergeometric (ordinary and $q$)
multisum/integral identities}, 
Invent. Math. 108 (1992) 575--633
\MR{1163239}

\bibitem[Wi]{Wi}
{\bf E Witten},
{\em Quantum field theory and the Jones polynomial}, 
Commun. Math. Physics. 121 (1989) 351--399
\MR{0990772}

\bibitem[Y]{Y}
{\bf T Yoshida},
{\em On ideal points of deformation curves of hyperbolic $3$--manifolds
with one cusp},
Topology 30 (1991) 155--170
\MR{1098911}

\bibitem[Z]{Z}
{\bf D Zeilberger},
{\em A holonomic systems approach to special functions identities},
J. Comput. Appl. Math. 32 (1990) 321--368
\MR{1090884}

\end{thebibliography}
\end{document}